\documentclass[12pt]{amsart}

\usepackage{graphicx}
\usepackage{amsmath}
\usepackage{amssymb}
\usepackage{amsfonts}
\usepackage{verbatim}
\usepackage{scalefnt}
\usepackage{multirow}
\usepackage{enumerate}
\usepackage{mathtools}
\usepackage{float}
\usepackage{enumitem}
\restylefloat{figure}
\usepackage[margin=3cm]{geometry}

\usepackage{fancyhdr}
\usepackage{hyperref}

\newcommand{\F}{{\mathbb F}}
\newcommand{\Z}{{\mathbb Z}}

\DeclareMathOperator{\slrk}{slice-rank}

\newtheorem{theorem}{Theorem}
\newtheorem{lemma}[theorem]{Lemma}

\newtheorem{corollary}[theorem]{Corollary}
\newtheorem{conjecture}[theorem]{Conjecture}

\newtheorem*{definition*}{Definition}

\newtheorem{remark}[theorem]{Remark}

\numberwithin{equation}{section}

\title[Avoiding right angles and certain Hamming distances]%
  {Avoiding right angles and certain Hamming distances}

\author{Bal\'azs Bursics}
\email{bursicsb@gmail.com}
\address{Eötvös Loránd University Faculty of Science, 1117 Budapest, Pázmány Péter sétány 1/A, Hungary}

\author{D\'avid Matolcsi}
\email{matolcsidavid@gmail.com}
\address{Eötvös Loránd University Faculty of Science, 1117 Budapest, Pázmány Péter sétány 1/A, Hungary}

\author{P\'eter P\'al Pach}
\email{ppp@cs.bme.hu}
\address{MTA-BME Lend\"ulet Arithmetic Combinatorics Research Group, Department of Computer Science and Information Theory, Budapest
  University of Technology and Economics, 1117 Budapest, Magyar tud\'osok
 k\"or\'utja 2., Hungary}

\author{Jakab Schrettner}
\email{sch.jakab@gmail.com}
\address{University of Cambridge, Centre for Mathematical Sciences, Wilberforce Rd, Cambridge CB3 0WA}

\thanks{}

\begin{document}

\begin{abstract}	
In this paper we show that the largest possible size of a subset of $\mathbb{F}_q^n$ avoiding right angles, that is, distinct vectors $x,y,z$ such that $x-z$ and $y-z$ are perpendicular to each other is at most $O(n^{q-2})$. This improves on the previously best known bound due to Naslund \cite{Naslund} and refutes a conjecture of Ge and Shangguan \cite{Ge}. A lower bound of $n^{q/3}$ is also presented.

It is also shown that a subset of $\mathbb{F}_q^n$ avoiding triangles with all right angles can have size at most $O(n^{2q-2})$. Furthermore, asymptotically tight bounds are given for the largest possible size of a subset $A\subseteq \mathbb{F}_q^n$ for which $x-y$ is not self-orthogonal for any distinct $x,y\in A$. The exact answer is determined for $q=3$ and $n\equiv 2\pmod {3}$.

Our methods can also be used to bound the maximum possible size of a binary code where no two codewords have Hamming distance divisible by a fixed prime $q$. Our lower- and upper bounds are asymptotically tight and both are sharp in infinitely many cases.

\end{abstract}

\date{\today}
\maketitle

\section{Introduction}

In this paper we consider problems about determining the largest possible size of sets avoiding certain geometric configurations. Some of our questions are related to coding theoretic problems.

Let $q$ be an odd prime power and $n$ a positive integer. A {\it right angle} in $\mathbb{F}_q^n$ is a triple  $x,y,z\in \mathbb{F}_q^n$ of distinct elements satisfying 
$$\langle x-z,y-z\rangle =0,$$
where $\langle \cdot,\cdot \rangle$ denotes the standard dot product.

Let $R(n,q)$ denote the largest possible size of a subset of $\mathbb{F}_q^n$ which contains no right angle. 

Bennett \cite{Bennett} proved that 
$$R(n,q)\ll q^{\frac{n+2}{3}}.$$
Ge and Shangguan \cite{Ge} used the so-called slice rank method to improve this bound for fixed $q$ and large $n$, namely, they showed that 
$$R(n,q)\leq \binom{n+q}{q-1}+3,$$
that is, for fixed $q$ the quantity $R(n,q)$ is only polynomial in $n$. 
They mentioned that the standard orthonormal basis yields the lower bound $R(n,q)\geq n$ and conjectured that their upper bound is asymptotically tight (for fixed $q$):
\begin{conjecture}[Ge-Shangguan, \cite{Ge}]\label{conj-GeS}
For any fixed prime power $q$, $R(n,q)=\Theta(n^{q-1})$.
\end{conjecture}
Naslund \cite{Naslund} further improved on the upper bound for $R(n,q)$ by showing that
$$R(n,q)\leq \binom{n+q}{q-1}+2-\binom{n+q}{q-3}.$$

We show that for every odd prime power $q$, in fact, $R(n,q)\ll n^{q-2}$, this refutes Conjecture~\ref{conj-GeS}:

\begin{theorem}\label{thm-right-angle}
Let $q$ be an odd prime power. If a set $A\subseteq \mathbb{F}_q^n$ does not contain a right angle, i.e. three distinct vectors with $\langle x-z,y-z\rangle=0$, then
$$|A|\leq 4(q-1)q\binom{n+q-2}{q-2}+2q.$$
\end{theorem}

The above problem is the finite
field version of the Erd\H{o}s-Falconer problem, which was originally defined in the setting of
Euclidean spaces and asked for the smallest $d$ for which any compact set in $\mathbb{R}^n$ with
Hausdorff dimension larger than $d$ contains three points forming an angle $\alpha$ (for given $n$ and $\alpha$). While in the real setting the question is interesting for any $\alpha$, in the case of the finite field version it is crucial that the dot product 0 is the forbidden one: If $A$ is an isotropic subspace of $\mathbb{F}_q^n$ of maximum dimension, then $\langle x-z,y-z\rangle=0$ for any $x,y,z\in A$, thus a set $A$ for which $\langle x-y,x-z\rangle\ne \alpha$ for some $\alpha\ne 0$ can be as large as $q^{\frac{n}{2}-1}$.

Naslund also considered a generalized version of right angles. Namely, we say that the vectors $x_0,x_1,\dots,x_k$ form a {\it $k$-right corner} if they are distinct, and if the $k$ vectors $x_1-x_0,\dots,x_k-x_0$ form a mutually orthogonal $k$-tuple, that is, $\langle x_i-x_0,x_j-x_0\rangle =0$ for all $1\leq i<j\leq k$. (Specially, a $2$-right corner is a right angle.) Naslund proved that for an integer $k$, an odd prime power $q=p^r$ with $p>k$ if a subset $A\subseteq \mathbb{F}_q^n$ satisfies
\begin{equation}\label{kright}
|A|>\binom{n+(k-1)q}{(k-1)(q-1)},
\end{equation}
then $A$ contains a $k$-right corner. 

We give the following lower bound for this problem:

\begin{theorem}\label{thm-Naslund-lower} Let $q$ be an odd prime and $2\leq k$ an integer. There exists a subset $A\subset \mathbb{F}_q^n$ of size
$$|A|\ge (1-o(1))\cdot \binom{n}{\big\lceil \frac{k-1}{k} \lfloor \frac{k}{2k-1}q \rfloor \big \rceil}\Big/ \binom{\lfloor \frac{k}{2k-1}q \rfloor}{\big\lceil \frac{k-1}{k} \lfloor \frac{k}{2k-1}q \rfloor \big \rceil}$$
which does not contain any $k$-right corner, i.e. vectors $x_0, x_1,\ldots, x_k$ such that 
$$\langle x_i-x_0, x_j-x_0 \rangle=0$$
for all $1\le i<j\le k$.
\end{theorem}

Note that for fixed $q$ and $k$ the upper bound for $k$-right corner free sets is of order $\Theta(n^{(k-1)(q-1)})$, while our lower bound is of order $\Theta\left(n^{\big\lceil \frac{(k-1)(q-1)}{2k-1}  \big \rceil}\right)$, where the exponent is the closest integer to $\frac{k-1}{2k-1}\cdot q$.

In the special case $k=2$ Theorem~\ref{thm-Naslund-lower} yields the following lower bound for $R(n,q)$:

\begin{corollary}\label{thm-ge-lower}
Let $q$ be an odd prime. There exists a subset $A\subset \mathbb{F}_q^n$ of size 
$$|A|\geq (1-o(1))\cdot \binom{n}{\big\lceil \frac{1}{2} \lfloor \frac{2}{3}q \rfloor \big \rceil}\Big/ \binom{\lfloor \frac{2}{3}q \rfloor}{\big\lceil \frac{1}{2} \lfloor \frac{2}{3}q \rfloor \big \rceil}$$
which does not contain a right angle, that is, it is not possible to choose distinct $x,y,z\in A$ such that $\langle y-x,z-x\rangle=0$.  
\end{corollary}
That is, $R(n,q)\gg n^{\big\lceil \frac{q-1}{3}  \big \rceil}$. Note that the exponent in this bound is the closest integer to $q/3$.

Motivated by these problems we also consider the following question: How large a set $A\subseteq \mathbb{F}_q^n$ can be if it contains no triangle with {\it all} right angles, that is, there are no distinct vectors $x,y,z\in A$ with
$$\langle x-y,y-z\rangle = \langle y-z, z-x \rangle = \langle z-x, x-y \rangle =0?$$

Using the slice-rank method, we obtain the following upper bound:

\begin{theorem}\label{thm-no-all-right}
Let $q$ be an odd prime. If $A\subseteq \mathbb{F}_q^n$  contains no triangle with \emph{all} right angles, i.e. vectors $x,y,z$ with 
$$\langle x-y,y-z\rangle = \langle y-z, z-x \rangle = \langle z-x, x-y \rangle =0,$$
then 
$$|A|\le \binom{n+2q-1}{2q-2}+2\binom{n+q}{q-1}.$$

\end{theorem}

It can be easily shown that the problem of avoiding triangles with all right angles is equivalent to avoiding triples $\{x,y,z\}$, where $$\langle x-y,x-y\rangle=\langle y-z,y-z\rangle =\langle z-x,z-x\rangle=0.$$
We consider the following related problem: How large can a set $A\subseteq \mathbb{F}_q^n$  be if $\langle x-y, x-y\rangle\ne0$ for any two distinct $x,y\in A$? For the maximal size of such a set in $\mathbb{F}_q^n$, we give the following bounds:

\begin{theorem}\label{thm-no-zero-length}

Let $q$ be an odd prime. 
Let $S(n,q)$ be the maximal size of a set in $\mathbb{F}_q^n$ which does not contain distinct vectors $x,y$ such that $\langle x-y, x-y\rangle=0$. We have the following bounds:
$$\binom{n}{q-1}\le S(n,q)\le \binom{n+q}{q-1}-\binom{n+q-2}{q-3}.$$
Moreover, whenever $n\not\equiv -2 \pmod q$ or $q \equiv 1 \pmod 4$, then 
$$S(n,q)\ge \binom{n}{q-1}+\binom{n}{q-2}.$$
\end{theorem}

In the special case $q=3$ we can further improve on the lower bound of Theorem~\ref{thm-no-zero-length} assuming that $n$ has residue 2 mod 3. Namely, we can construct a suitable set of size $\binom{n+3}{2}-1$ which is our general upper bound. Therefore, in infinitely many cases the exact answer is determined according to the following:

\begin{theorem}\label{thm-mod3-exact-value}
For $n\equiv 2 \pmod 3$ we have 
\[S(n,3)= \binom{n+3}{2}-1.\]
\end{theorem}

For $n\equiv 0 \pmod 3$ resp. $n\equiv 1 \pmod 3$ we can achieve the lower bounds $\binom{n+2}{2}-1$ resp. $\binom{n+1}{2}-1$ by choosing the last 1 resp. 2 digits to be constant and taking the construction from Theorem~\ref{thm-mod3-exact-value} on the remaining entries. For $n=3,4$ these turn out to be sharp, however, for infinitely many values of $n$ there is a better construction. Namely, whenever $q\equiv 2 \pmod 3$ is a prime power, and $n=q^2+q+1$, we have a construction of size $\binom{n+1}{2}$, which is bigger (by one!) than the construction given above.

We can also think of the problem of determining $S(n,3)$ from a coding theoretic point of view. Since in $\mathbb{F}_3$ the square of every nonzero element is 1, the value of $\langle x-y,x-y\rangle$ is the same as the Hamming-distance of $x$ and $y$ (modulo 3). Therefore, $S(n,3)$ is the largest possible size of a ternary code where none of the Hamming-distances between the codewords is divisible by 3.

In the construction from the proof of Theorem~\ref{thm-mod3-exact-value} we use the following observation: There are $\binom{n}{2}+1$ vectors in $\{0,1\}^n(\subseteq \mathbb{F}_3^n)$ with no two of them having Hamming distance divisible by 3. We investigate this question more generally: What is the largest possible size of a binary code of length $n$, if the Hamming distance of two different codewords is never divisible by a fixed prime $q$?

We provide the following bounds:

\begin{theorem}\label{thm-Hamming}
Let $q$ be an odd prime.
Let $T(n,q)$ be the maximal size of a subset of $\{a,b\}^n$ with no two vectors having Hamming distance divisible by a fixed prime $q$. Then the following bounds hold:
\begin{align*}
T(n,q)&\le \binom{n}{q-1}+\binom{n}{q-2}+\dots+\binom{n}{1}+\binom{n}{0}\qquad \text{in general,}\\
T(n,q)&\le \binom{n}{q-1}+\binom{n}{q-3}+\dots+\binom{n}{2}+\binom{n}{0} \qquad\text{for } n\equiv 0 \pmod q,\\
T(n,q)&\ge\binom{n}{q-1}+\binom{n}{q-3}+\dots+\binom{n}{2}+\binom{n}{0}\qquad \text{in general,}\\
T(n,q)&\ge\binom{n}{q-1}+\binom{n}{q-2}+\dots+\binom{n}{1}+\binom{n}{0}\qquad\text{for } n\equiv -1 \pmod q.
\end{align*}
\end{theorem}

Hence, by Theorem~\ref{thm-Hamming} the exact values are determined for $n\equiv 0\text{ or }-1 \pmod{q}$. 

For the particular case $q=3$, Theorem~\ref{thm-Hamming} leaves open only the case $n\equiv 1 \pmod 3$. In this case, for $n=4$ the exact value is 8, which is 1 larger than the general lower bound, but we do not know the exact values for any larger $n$ which has residue 1 mod 3. 

Note that the general upper bound from Theorem~\ref{thm-Hamming} is an old result of Delsarte \cite{Delsarte}. Different proofs were given by Frankl \cite[Theorem 1.6]{Frankl} and Babai et el \cite{Babai}. In this paper the prime case is considered, we give yet another proof of the statement, and provide a general lower bound. Furthermore, we show that both the lower and upper bounds are tight in infinitely many cases.


\section{Notation and preliminaries}

\subsection{Notation.}
Throughout the paper, the standard notation $\ll$, $\gg$
and respectively $O$ and $\Omega$ is applied to positive quantities in the usual way.
That is, $X \gg Y$, $ Y \ll X$, $X = \Omega(Y )$ and $Y = O(X)$ all mean that $X \geq cY$,
for some absolute constant $c > 0$. If both $X \ll Y$ and $Y \ll X $ hold, we write
 $X = \Theta(Y )$. If the constant $c$ depends on a quantity
$t$, we write $X \ll_t Y$, $Y = O_t(Y )$, and so on.

Throughout the paper the standard dot product of vectors $x,y\in \mathbb{F}_q^n$ is denoted by $\langle x,y\rangle=\sum\limits_{i=1}^n x_iy_i$. 

We will use the notation $[n]:=\{1,2,\dots,n\}$ and the Kronecker delta function as
$$\delta_{x,y}=\delta_x(y)=\delta(x,y)=\begin{cases} 1 & \text{if }x=y\\
0 & \text{if }x\ne y\end{cases}.$$

\subsection{Slice rank method.}

For proving some of the upper bounds we shall use the so-called {\it slice rank method}. This is the symmetrized version given by Tao of the polynomial method developed by Croot, Lev and the third named author \cite{Croot-Lev-Pach} to show that 3AP-free sets in $\mathbb{Z}_4^n$ are exponentially small which was subsequently adapted by Ellenberg and Gijswijt \cite{Ellenberg-Gijswijt} to the case of $\mathbb{F}_3^n$. For an exposition of the slice rank method we refer to \cite{Tao1,Tao2}. Here, let us briefly summarize the method. For finite sets $A_1,\dots,A_k$ and a field $\mathbb{F}$ we say that
$$h:A_1\times \dots \times A_k\to \mathbb{F}$$
has slice-rank 1 if 
$$h(x_1,\dots,x_k)=f(x_i)g(x_1,\dots,x_{i-1},x_{i+1},\dots,x_k)$$
for some $1\leq i\leq k$. The slice rank of $F:A_1\times \dots \times A_k\to \mathbb{F}$, denoted $\slrk (F)$, is the smallest $r$ such that $F=\sum\limits_{i=1}^r h_i$ where the $h_i$ functions have slice rank 1.

The following lemma from \cite{Tao1} is going to be applied several times:
\begin{lemma}\label{lem-diag-slice}
Let $A$ be a finite set and let $F:\underbrace{A\times \dots \times A}_k\to \mathbb{F}$ be a diagonal tensor attaining only nonzero entries on its diagonal, that is, 
$$F(x_1,\dots,x_k)=\sum\limits_{a\in A} c_a \delta_a(x_1)\dots \delta_a(x_k)$$
where $c_a\ne 0\ (a\in A)$ and $\delta_a(x)=\begin{cases} 1 & x=a\\
0 & x\ne a\end{cases}$.

Then 
$$\slrk (F)=|A|.$$
\end{lemma}
In the applications for bounding the size of a set $A$ satisfying certain criteria, a diagonal tensor $F$ is going to be presented, then it will suffice to bound the slice-rank of $F$.

For bounding the slice-rank the tensor $F$ is going to be represented as a multivariate polynomial, and in certain cases we will have to calculate with the dimension of certain subspaces of polynomials. Namely, we will use multiple times that the number of monomials in $\mathbb{F}[x_1,\dots,x_n]$ with total degree at most $d$ is $\binom{n+d}{d}$. Indeed, the number of $n$-variable monomials with total degree exactly $i$ is $\binom{n+i}{i-1}$, and $\sum\limits_{i=0}^d \binom{n+i}{i-1}=\binom{n+d}{d}$. 

One application of the slice rank method is going to be through the following lemma:

\begin{lemma}\label{lem-few-dot-product}
Let $q$ be a prime power. 
Let $\alpha\in \mathbb{F}_q$ and $R\subseteq \mathbb{F}_q\setminus \{\alpha\}$. Let us assume that for a set $A\subseteq \mathbb{F}_q^n$  we have $\langle a,a\rangle =\alpha$ for every $a\in A$ and $\langle x,y\rangle\in R$ for any two distinct $x,y\in A$. Then $|A|\leq 2\binom{n+|R|}{|R|}$. 
\end{lemma}

\begin{proof}
Let us consider $T:\mathbb{F}_q^n\times \mathbb{F}_q^n \to \mathbb{F}_q$ defined as 
$$T(x,y)=\prod\limits_{r\in R} (\langle x,y\rangle -r).$$
Note that $T(x,x)\ne 0$ and $T(x,y)=0$ for any $x\ne y$, $x,y\in R$. By Lemma~\ref{lem-diag-slice} we have $|A|=\slrk T\big|_{A\times A}\leq \slrk T$.

As $\deg T=2|R|$, the slice-rank of $T$ is at most 2 times the number of those monomials with $n$ variables that have degree at most $|R|$. Indeed, in each monomial either the total degree of the $x$'s or the total degree of the $y$'s is at most $|R|$. Thus $T$ may be expressed as a sum of products $fg$, where
\begin{itemize}
    \item either $f$ is a monomial of total degree at most $|R|$ from $\mathbb{F}_q[x_1,x_2,\dots,x_n]$  and $g\in \mathbb{F}_q[y_1,y_2,\dots,y_n]$,
    \item or $f$ is a monomial of total degree at most $|R|$ from $\mathbb{F}_q[y_1,y_2,\dots,y_n]$  and $g\in \mathbb{F}_q[x_1,x_2,\dots,x_n]$.
\end{itemize}

As the number of $n$-variable monomials of total degree at most $|R|$ is $\binom{n+|R|}{|R|}$, we obtain that
$$|A|\leq \slrk T\leq 2\binom{n+|R|}{|R|}.$$
\end{proof}

\subsection{An extremal set theoretic lemma.}

For our lower bound constructions we shall use the following lemma about $t$-uniform families of sets with bounded intersection size:

\begin{lemma}\label{lem-t-l-int}
Let $\ell\leq t$ be positive integers. 
There exists a $t$-uniform family $\mathcal A$ of subsets of $[n]$ such that for any two distinct $F,G\in \mathcal A$ we have $|F\cap G|<\ell$ and 
$$|\mathcal A|\geq (1-o_{n\to\infty}(1))\cdot\frac{\binom{n}{\ell}}{\binom{t}{\ell}}.$$
\end{lemma}

\begin{proof}
The statement was proved by R\"odl  in \cite{Rodl}.
\end{proof}

Note that to achieve a lower bound that is valid for every $n$ one can generalize \cite[Proposition 1]{EFF} to get a construction of size $|\mathcal A|\geq \frac{\binom{n}{\ell}}{\binom{t}{\ell}^2}.$

\section{Proofs}


\begin{proof}[Proof of Theorem~\ref{thm-right-angle}]
Let $A\subseteq \mathbb{F}_q^n$ be a set such that $\langle x-z,y-z\rangle \ne 0$ for distinct $x,y,z\in A$. For some $\alpha\in \mathbb{F}_q$ for at least $|A|/q$ elements of $A$ we have $\langle x,x\rangle =\alpha$. Let $A^{\alpha}:=\{x\in A:\langle x,x\rangle=\alpha\}$.

Let us pick an element $u\in A^{\alpha}$. For $\beta\in \mathbb{F}_q$ let $$A^{\alpha}_{\beta}=\{x\in A^{\alpha}\setminus\{u\}: \langle u,x\rangle =\beta \}.$$
For any two distinct $x,y\in A^{\alpha}_{\beta} $ we have $\langle x,y\rangle \in \mathbb{F}_q\setminus \{ \alpha,2\beta-\alpha\}$. Indeed, $\langle x,y\rangle =\alpha$ would imply $\langle x-y,u-y\rangle =0$, while $\langle x,y\rangle =2\beta-\alpha$ would imply $\langle x-u,y-u\rangle =0$.

For $\beta\ne \alpha$ Lemma~\ref{lem-few-dot-product} applied with the choice $R:=\mathbb{F}_q\setminus \{\alpha,2\beta-\alpha\}$ readily implies that 
$$|A^{\alpha}_{\beta}|\leq 2\binom{n+q-2}{q-2}.$$ 

To also bound $A^{\alpha}_\alpha$, let us
pick an element $v\in A^{\alpha}_{\alpha}$. (Note that $v\ne u$.) For   $\beta\in \mathbb{F}_q$ let 
$$A^{\alpha}_{\alpha,\beta}=\{x\in A^{\alpha}_{\alpha}\setminus\{v\}: \langle v,x\rangle =\beta \}.$$
Note that $A^{\alpha}_{\alpha,\alpha}=\emptyset$ and again by Lemma~\ref{lem-few-dot-product} 
$$|A^{\alpha}_{\alpha,\beta}|\leq 2\binom{n+q-2}{q-2}$$
for every $\beta\ne \alpha$. 

Therefore, 
$$|A^{\alpha}|=1+|A^{\alpha}_\alpha|+\sum\limits_{\beta\in \mathbb{F}_q\setminus \{\alpha\}} |A^{\alpha}_\beta|=2+\sum\limits_{\beta\in \mathbb{F}_q\setminus \{\alpha\}} |A^{\alpha}_{\alpha,\beta}|+\sum\limits_{\beta\in \mathbb{F}_q\setminus \{\alpha\}} |A^{\alpha}_\beta|\leq 4(q-1)\binom{n+q-2}{q-2}+2.$$
Hence, $|A|\leq q|A^{\alpha}|\leq 4(q-1)q\binom{n+q-2}{q-2}+2q$.
\end{proof}

\begin{remark}
Note that we did not optimize in the proof the arising constant factor $4(q-1)q$. 

 
\end{remark}

\begin{proof}[Proof of Theorem~\ref{thm-Naslund-lower}]

Let $t=\lfloor \frac{k}{2k-1}q \rfloor$. Let $\mathcal{A}$ be a $t$-uniform family of subsets of $[n]$ such that for any two distinct $F,G \in \mathcal{A}$ we have $|F\cap G|<\frac{k-1}{k}t$. According to Lemma~\ref{lem-t-l-int} we may take a system of size $|\mathcal A|\geq (1-o(1))\cdot \frac{\binom{n}{\left\lceil (k-1)t/k \right\rceil}}{\binom{t}{\left\lceil (k-1)t/k \right\rceil}}$.

Now, let $A$ consist of the characteristic vectors of the elements of $\mathcal A$, considered as elements of $\mathbb{F}_q^n$. That is, $a=(a_1,\dots,a_n)\in \{0,1\}^n$ is contained in $A$ if and only if $\{i:\ a_i=1\}\in \mathcal A$.

We claim that $A$ avoids $k$-right corners, that is, for any $k+1$ distinct elements $x_0,x_1,\dots,x_k\in A$ the dot products $\langle x_i-x_0,x_j-x_0\rangle$ ($1\leq i<j\leq k$) can not be simultaneously 0 (in $\mathbb{F}_q$).

To show this, let us take arbitrarily $k+1$ distinct elements, $x_0,x_1,\dots,x_k\in A$, and denote by $X_0,X_1,\dots,X_k$ the corresponding subsets of $[n]$. First, let us show that for some indices $1\leq i<j\leq k$ we have $X_0\setminus (X_i\cup X_j)\ne \emptyset$. Since $X_0\in \mathcal A$, we have $|X_0|=t$, and we also know that $|X_0\setminus X_i|>\frac{t}{k}$ for every $1\leq i\leq k$. Thus, by the pigeonhole principle we obtain that for two different indices $i$ and $j$ we have $\emptyset \ne (X_0\setminus X_i)\cap (X_0\setminus X_j)=X_0\setminus (X_i \cup X_j)$, as we claimed.

Observe that
\begin{multline*}\langle x_i-x_0,x_j-x_0\rangle =\langle x_i,x_j\rangle +\langle x_0,x_0\rangle -\langle x_i,x_0\rangle -\langle x_0,x_j\rangle=\\
=|X_i\cap X_j|+|X_0|-|X_i\cap X_0|-|X_0\cap X_j|=|(X_i\cap X_j)\setminus X_0|+|X_0\setminus (X_i\cup X_j)|.
\end{multline*}
As $0\leq |(X_i\cap X_j)\setminus X_0|<\frac{k-1}{k}t$ and $0<|X_0\setminus (X_i\cup X_j)|\leq t$ we have
$$0<|(X_i\cap X_j)\setminus X_0|+|X_0\setminus (X_i\cup X_j)|<\frac{k-1}{k}t+t=\frac{2k-1}{k}t\leq q,$$

hence $\langle x_i-x_0,x_j-x_0\rangle \ne 0$.

\end{proof}



\begin{proof}[Proof of Theorem~\ref{thm-no-all-right}]

First of all, we show that the condition that the three sides of a triangle are pairwise orthogonal is equivalent to the condition that all the three sides are self-orthogonal. Assume first that 
$$\langle x-y,y-z\rangle = \langle y-z, z-x \rangle = \langle z-x, x-y \rangle =0.$$ 
Now, 
$$\langle x-y,x-y\rangle = \langle x-y,(x-z)+(z-y)\rangle=\langle x-y,x-z\rangle+\langle x-y,z-y\rangle=0+0=0,$$
and $\langle y-z,y-z\rangle=\langle z-x,z-x\rangle=0$ can be proved analogously.

For the reverse direction assume that 
$$\langle x-y,x-y\rangle=\langle y-z,y-z\rangle=\langle z-x,z-x\rangle=0.$$
Then
$$\langle x-y,y-z\rangle =\frac{ \langle (x-y)+(y-z),(x-y)+(y-z)  \rangle -\langle x-y,x-y\rangle -\langle y-z,y-z\rangle }{2}=0,$$
and $\langle x-y,z-x\rangle =\langle y-z,z-x\rangle=0$ follows similarly. 

Hence, if $A\subseteq \mathbb{F}_q^n$ does not contain a triangle with all right angles, then the set $A$ does not contain three distinct vectors $x,y,z$ such that $\langle x-y,x-y\rangle =\langle y-z,y-z\rangle= \langle z-x,z-x\rangle=0$.

Let us consider  $F:\mathbb{F}_q^n\times\mathbb{F}_q^n\times\mathbb{F}_q^n\to \mathbb{F}_q$ defined as
$$F(x,y,z)=(1-\langle x-y,x-y\rangle^{q-1})(1-\langle y-z,y-z\rangle^{q-1})(1-\langle z-x,z-x\rangle^{q-1}).$$
Observe that for distinct $x,y,z\in A$ we have $F(x,y,z)=0$ according to the assumption on $A$. Also, for $x=y=z$ we have $F(x,y,z)=F(x,x,x)=1$. However, $F\big|_{A\times A\times A}$ might not be diagonal, as for instance $F(x,x,z)$ may be nonzero for $x\ne z.$ Therefore, let us consider $$G(x,y,z):=F(x,y,z)(1-\delta(x,y)-\delta(y,z)-\delta(z,x)),$$ where $\delta$ is the Kronecker delta function. For $x=y=z$ we have $G(x,y,z)=1\cdot (1-1-1-1)=-2\ne 0$ and $G$  vanishes at $(x,y,z)$ when $x,y,z$ are distinct elements of $A$. Moreover, if two of $x,y,z$ coincide, but the third one is different, then $1-\delta(x,y)-\delta(y,z)-\delta(z,x)=0$, thus $G(x,y,z)=0$ also holds.

Hence, $G\big|_{A\times A\times A}$ is diagonal, that is, for $x,y,z\in A$ we have $G(x,y,z)\ne 0$ if and only if $x=y=z$. According to Lemma~\ref{lem-diag-slice} we have $|A|=\text{slice-rank}(G\big|_{A\times A\times A})$. Now, we give an upper bound for the slice-rank of $G\big|_{A\times A\times A}$. 

First of all,
\begin{multline*}F(x,y,z)=(1-\langle x-y,x-y\rangle^{q-1})(1-\langle y-z,y-z\rangle^{q-1})(1-\langle z-x,z-x\rangle^{q-1})=\\
=[1-(x_0+y_0-2x_1y_1-2x_2y_2-\dots-2x_ny_n)^{q-1}]\cdot[1-(y_0+z_0-2y_1z_1-2y_2z_2-\dots-2y_nz_n)^{q-1}]\cdot\\
\cdot[1-(z_0+x_0-2z_1x_1-2z_2x_2-\dots-2z_nx_n)^{q-1}],
\end{multline*}
where $x_0=\sum\limits_{i=1}^n x_i^2$, $y_0=\sum\limits_{i=1}^n y_i^2$ and $z_0=\sum\limits_{i=1}^n z_i^2$. 

Let us consider $x_0$ as a new variable. Each monomial in $F(x,y,z)$ can be written as a product of a monomial in $\mathbb{F}_q[x_0,x_1,\dots,x_n]$ of degree at most $2(q-1)$ and a polynomial of $y$ and $z$. Therefore, the slice-rank of $F$ is at most the number of $(n+1)$-variable monomials of degree at most $2(q-1)$, which is $\binom{n+2q-1}{2q-2}$.

Let us consider now $\delta(x,y)F(x,y,z)$. 
Observe that on $A\times A\times A$ we have
$$\delta(x,y)F(x,y,z)=\delta(x,y)(1-\langle z-x,z-x\rangle^{q-1}).$$
Note that each monomial in $1-\langle z-x,z-x\rangle^{q-1}=1-(z_0+x_0-2z_1x_1-\dots-2z_nx_n)^{q-1}$ can be written as a product of a monomial in $\mathbb{F}_q[z_0,z_1,\dots,z_n]$ of degree at most $q-1$ and a polynomial of $x$ and $y$. The number of $(n+1)$-variable monomials of degree at most $q-1$  is $\binom{n+q}{q-1}$. Since $\delta(x,y)=\prod\limits_{i=1}^n (1-(x_i-y_i)^{q-1})$ is a polynomial of $x$ and $y$, the slice-rank of $\delta(x,y)F(x,y,z)\big|_{A\times A\times A}$ is at most $\binom{n+q}{q-1}$.

Analogously, the slice-rank of $\delta(y,z)F(x,y,z)\big|_{A\times A\times A}$  is also at most $\binom{n+q}{q-1}$. Also, each monomial of $\delta(z,x)F(x,y,z)\big|_{A\times A\times A}$ can be written as a product of a mo\-no\-mial in $\mathbb{F}_q[x_0,x_1,\dots,x_n]$ of degree at most $q-1$ and a polynomial of $y$ and $z$. However, these types of products already appeared when we considered the slice-rank $F$, so in fact 
$$\slrk (F(x,y,z)(1-\delta(x,y))\leq \binom{n+2q-1}{2q-2}.$$

Hence, 
$$|A|= \text{slice-rank}\left(G\big|_{A\times A\times A}\right)\leq \binom{n+2q-1}{2q-2}+2\binom{n+q}{q-1}.$$

\end{proof}




\begin{proof}[Proof of Theorem~\ref{thm-no-zero-length}] 

Let us start with proving the upper bound.

Let $A$ be a set satisfying the property, that is, for every pair of distinct $x,y\in A$ we have $\langle x-y,x-y\rangle\ne 0$. For every $a\in A$ let $p_a$ be the following polynomial:  
$$p_a(x):=1-\langle x-a,x-a \rangle ^{q-1}=1-\left(\sum\limits_{i=1}^n (x_i-a_i)^2 \right)^{q-1}.$$

Note that $p_a(a)=1$, and  for every $b\in A\setminus \{a\}$ we have $p_a(b)=0$, since $\langle a-b,a-b\rangle\ne 0$. This readily implies that the system of polynomials $\{p_a\ |\ a\in A\}$ is linearly independent.

Now, we prove that the polynomials $p_a$ are all contained in a subspace of dimension $\binom{n+q}{q-1}-\binom{n+q-2}{q-3}$. In order to see this let us define the following sets of polynomials:

$$H=\Big\{ (x_1^2+ \dots +x_n^2)^{\alpha_0}\cdot x_1^{\alpha_1} \cdot x_2^{\alpha_2}\cdot \ldots \cdot x_n^{\alpha_n}\ |\ \displaystyle\sum_{i=0}^n \alpha_i \le q-1 \Big\},$$

$$I=\Big\{x_1^{\alpha_1} \cdot x_2^{\alpha_2}\cdot \ldots \cdot x_n^{\alpha_n}\ |\ \displaystyle\sum_{i=1}^n \alpha_i \le q-1 \Big\},$$

$$J=\Big\{ (x_1^2+ \dots +x_n^2)^{\alpha_0}\cdot x_1^{\alpha_1} \cdot x_2^{\alpha_2}\cdot \ldots \cdot x_n^{\alpha_n}\ |\ \alpha_0 \ge 1, \displaystyle\sum_{i=0}^n \alpha_i=q-1 \Big\}.$$

Since 
\begin{multline*}
p_a(x)=1-\left(\sum\limits_{i=1}^n (x_i-a_i)^2 \right)^{q-1}=\\
=1-\left((x_1^2+\dots+x_n^2)-2a_1x_1-\dots-2a_nx_n+(a_1^2+\dots+a_n^2) \right)^{q-1},
\end{multline*}
it is clear that
 $\{p_a\ |\ a\in A\} \subseteq \mathrm{span}_{\mathbb{F}_q}(H)$. We will show that  $\mathrm{span}_{\mathbb{F}_q}(H) \subseteq \mathrm{span}_{\mathbb{F}_q}(I \cup J)$. 

Let $p=(x_1^2+ \dots +x_n^2)^{\alpha_0}\cdot x_1^{\alpha_1} \cdot x_2^{\alpha_2}\cdot \ldots \cdot x_n^{\alpha_n}\in H$. For the sake of contradiction, assume that $p\notin \mathrm{span}_{\mathbb{F}_q}(I \cup J)$. We may further assume that such a $p$ is chosen in such a way that $\alpha_0$ is minimal. Note that since $p\notin I\cup J$, we have $\alpha_0\geq 1$ and $\sum\limits_{i=0}^n \alpha_i<q-1$.
By the assumption on the minimality of $\alpha_0$ for every $1\leq i\leq n$ the polynomial 
$$p_i:=(x_1^2+\dots+x_n^2)^{\alpha_0-1}x_1^{\alpha_1}\dots x_i^{\alpha_i+2}\dots x_n^{\alpha_n}\in H$$
is contained in $\mathrm{span}_{\mathbb{F}_q}(I \cup J)$. However, this would imply that $p=p_1+\dots+p_n\in \mathrm{span}_{\mathbb{F}_q}(I \cup J)$, which contradicts our assumption.

Therefore, $\{p_a\ |\ a\in A\} \subseteq \mathrm{span}_{\mathbb{F}_q}(H) \subseteq \mathrm{span}_{\mathbb{F}_q}(I \cup J)$. Hence,
$$|A| \le \mathrm{dim}\big( \mathrm{span}_{\mathbb{F}_q}(H)\big) \le \mathrm{dim}\big( \mathrm{span}_{\mathbb{F}_q}(I \cup J)\big) \le |I|+|J|.$$

Since $|I|=\binom{n+q-1}{q-1}$ and $|J|=\sum\limits_{i=1}^{q-1} \binom{n+q-2-i}{q-1-i}=\binom{n+q-2}{q-2}$ we obtain that $|A|\leq \binom{n+q-1}{q-1}+\binom{n+q-2}{q-2}=\binom{n+q}{q-1}-\binom{n+q-2}{q-3}$, which completes the proof of the upper bound.


\medskip

Let us continue with the lower bound. 


	Let $A$ be the set of all vectors having $q-1$ components equal to 1, and all the remaining $n-q+1$ components being equal to 0. Note that $|A|=\binom{n}{q-1}$. 
	
	If $x,y\in A$ are distinct, then $\langle x-y,x-y\rangle =2(q-1)-2\langle x,y\rangle\ne 0$, since $\langle x,y\rangle\ne q-1$.

\medskip


Finally, in the cases when $n\not\equiv -2 \pmod q$ or $q \equiv 1 \pmod 4$ we can further improve on the previous lower bound. The construction is as follows: let $A_1$ be the previously defined set of all vectors with $q-1$ components being equal to 1, and $n-q+1$ components being equal to 0. Furthermore, let $A_2$ be the set of all vectors with $q-2$ entries $a$ and $n-q+2$ entries $b$, where the elements $a\not=b$ are going to be chosen later. Note that the size of $A$ is  $\binom{n}{q-1}+\binom{n}{q-2}$. We claim that for an appropriate choice of $a$ and $b$ the set $A=A_1\cup A_2$ does not contain two vectors $x,y$ such that $\langle x-y,x-y\rangle=0$.

We have already seen that $\langle x-y,x-y\rangle\ne 0$ if $x,y$ are distinct elements taken from $A_1$. Similarly, if $x,y$ are distinct vectors taken from $A_2$, then $\langle x-y,x-y\rangle =2k(a-b)^2\ne 0$, where $k$ is the number of entries that are $a$ in $x$ and $b$ in $y$. (As $a\ne b$, we have $0<k\leq q-2$, and by the definition of $A_2$ the number of those entries that are $b$ in $x$ and $a$ in $y$ is also $k$, all the remaining entries are the same in $x$ and $y$.)

It remains to consider the pairs $x\in A_1,\ y\in A_2$.

 Let us assume that there are $k$ indices $i$ with $(x_i,y_i)=(1,a)$. Then $0\le k\le q-2$ and there are
	\begin{itemize}
		\item $q-1-k$ indices $i$ with $(x_i,y_i)=(1,b)$,
		\item $q-2-k$ indices $i$ with $(x_i,y_i)=(0,a)$,
		\item $n+k-2q+3$ indices $i$ with $(x_i,y_i)=(0,b).$
	\end{itemize}
	Therefore, (calculating in $\mathbb{F}_q)$
	\begin{multline*}
		\langle x-y,x-y\rangle =k(a-1)^2+(q-1-k)(b-1)^2+(q-2-k)a^2+(n+k-2q+3)b^2=\\
		=-2a^2+(n+2)b^2-2ka+(2k+2)b-1=(2b-2a)k-2a^2+(n+2)b^2+2b-1.
	\end{multline*}

	Thus we can think of $\langle x-y,x-y\rangle$ as a linear function of $k$ with leading coefficient $-2a+2b=2(b-a)\not=0$ (since  $a\ne b$). Hence, there is exactly one residue class for ($k$ mod $q$) which solves this linear equation (and gives $\langle x-y,x-y\rangle=0$). Note that $k\in [0,q-2]$, so to ensure $\langle x-y,x-y\rangle\not=0$ for all choices of $x,y$ this (unique) solution must be  $k\equiv -1 \pmod q$. That is, we shall choose $a$ and $b$ in such a way that
	$$-2a^2+(n+2)b^2+2a-1=0,$$
	or equivalently
	\begin{equation}\label{eq-ab}
	(n+2)b^2=2a^2-2a+1.\end{equation}
	
Here we distinguish two cases depending on whether $n\not\equiv -2 \pmod q$ or $n\equiv -2 \pmod q$. Let us assume first that $n\not\equiv -2 \pmod q$.

Observe that both sides of the equation \eqref{eq-ab} can attain exactly $\frac{q+1}{2}$ different values, thus, there must be at least one element which can be expressed as $(n+2)b^2$ and as $2a^2-2a+1$, simultaneously. In other words, the equation has a solution $a,b$. If $a\ne b$, then we are done, otherwise we may take the solution $a,-b$, which also satisfies $a\ne -b$ (since $a=b=0$ is not a solution).

Finally, we consider the case $n\equiv -2 \pmod{q}$. Here, the solvability of \eqref{eq-ab} reduces to the solvability of $0=2a^2-2a+1$. This equation is solvable if its discriminant $-4$ is a square modulo $q$, that is, if $q\equiv 1 \pmod{4}$. In this case for a solution $a$ we may choose any $b\ne a$. This completes our proof.

\end{proof}


\begin{proof}[Proof of Theorem~\ref{thm-mod3-exact-value}]
The upper bound follows from Theorem~\ref{thm-no-zero-length}, by presenting a matching lower bound the proof will  be complete. The construction for the lower bound is as follows:

\noindent
 Let the set $A\subseteq \mathbb{F}_3^n$ consist of

\begin{itemize}
    \item all the vectors where exactly two entries are $1$'s and all other entries are $0$'s.
    \item all the vectors where one entry is 0 and all other entries are $2$'s,
    \item all the vectors where one entry is $1$ and all the other entries are $2$'s,
    \item all the vectors where one entry is $0$ and all the other entries are $1$'s,
    \item the all-zero vector $(0,0,\dots,0)$,
    \item the all-2 vector $(2,2,\dots,2)$.
\end{itemize}

Since we are modulo 3, we have $\langle x-y,x-y\rangle =0$ if and only if the Hamming-distance of the vectors $x$ and $y$ is divisible by 3. It is easy to check that assuming $n\equiv 2\pmod{3}$ none of the pairs of vectors from the above defined set have Hamming-distance divisible by 3.

The cardinality of the set is $|A|=\binom{n}{2}+3n+2=\binom{n+3}{2}-1$.
\end{proof}


\begin{proof}[Proof of Theorem~\ref{thm-Hamming}]

Let us start with proving the lower bound. Let $A$ consist of those elements of $\{a,b\}^n$ in which the number of characters $a$ is even and at most $q-1$.

Let $x$ and $y$ be two distinct elements of $A$. Let $k$ resp. $\ell$ denote the number of characters $a$ in $x$ resp. $y$, and denote by $m$ the number of those entries where both $x$ and $y$ have $a$. Then the number of indices $i$ such that $(x_i,y_i)=(a,b)$ is $k-m$ and the number of indices $i$ such that $(x_i,y_i)=(b,a)$ is $\ell-m$. Hence, the total number of those indices where $x$ and $y$ differ from each other is $k+\ell-2m$:
$$d(x,y)=k+\ell-2m,$$
where $d(x,y)$ stands for the Hamming distance of $x$ and $y$. 
Observe that $k$ and $\ell$ are even and $k,\ell\leq q-1$, thus $d(x,y)=k+\ell-2m\in (0,2q)$ is even, which implies that $d(x,y)$ is not divisible by $q$.

As $|A|= \binom{n}{q-1}+\binom{n}{q-3}+\dots+\binom{n}{2}+\binom{n}{0}$, this completes the proof of our general lower bound.

In the special case when $n\equiv -1 \pmod{q}$, we present a better construction.

Let us add to the previously defined set $A$ those elements of $\{a,b\}^n$ in which the number of characters $b$ is odd and at most $q-2$. This way the set $A'$ is obtained which has size $|A'|=\binom{n}{q-1}+\binom{n}{q-2}+\dots+\binom{n}{1}+\binom{n}{0}$.

We have already seen that $q\nmid d(x,y)=k+\ell-2m$ if $x,y\in A$ are distinct.

If $x$ and $y$ are two distinct elements of $A'\setminus A$, then $n-k$ and $n-\ell$ are odd and at most $q-2$. The number of those indices for which $(x_i,y_i)=(b,b)$ is $n-k-\ell+m$, thus 
$$0<d(x,y)=k+\ell-2m=(n-k)+(n-\ell)-2(n+m-k-\ell)\leq 2(q-2)<2q,$$
and $k+\ell-2m$ is even, therefore, $d(x,y)=k+\ell-2m\ne q$, so $q\nmid d(x,y)$.

Finally, let us assume that $x\in A$ and $y\in A'\setminus A$. Then $k\leq q-1$ is even and $n-\ell\leq q-2$ is odd. Observe that $$d(x,y)=k+\ell-2m=n+2(k-m)-k-(n-\ell),$$
thus 
$$n-2q+3=n+0-(q-1)-(q-2)\leq n+2(k-m)-k-(n-\ell)= d(x,y)\leq n.$$

Since $2\mid k$ and $2\nmid n-\ell$, the Hamming-distance $d(x,y)=k+\ell-2m$ has the same parity as $n+1$. As $q\mid n+1$, in the interval $[n-2q+3,n]$ only one element, $n-q+1$ is divisible by $q$, however, its parity is different from the parity of $n+1$. Hence, $q\nmid d(x,y)$, as we claimed.

\bigskip

Let us continue with the upper bound. We may assume that the two characters are $\pm1$, considered as elements of $\mathbb{F}_q$. This way, $A\subseteq \{-1,1\}^n\subset \mathbb{F}_q^n$. Note that for $x,y\in \{-1,1\}^n$ we have $\langle x-y,x-y\rangle =\sum\limits_{i=1}^n (x_i-y_i)^2=4d(x,y)$, as $(x_i-y_i)^2$ is 4 if $x$ and $y$ differ in the $i$th coordinate and 0 otherwise. Therefore, the Hamming-distance $d(x,y)$ is divisible by $q$ if and only if $\langle x-y,x-y\rangle=0$.

Let $A=\{a_1,a_2,\dots,a_r\}$. The assumption on $A$ implies that $\langle a_i-a_j,a_i-a_j\rangle\ne0$, when $i\ne j$. Let us consider the polynomials $p_i(x):=1-\langle x-a_i,x-a_i\rangle^{q-1}$. Let $f_i:=p_i\big|_{A}$, that is, by restricting the domain of $p_i$ to $A$ we obtain the function $f_i$. 
Note that for every $i,j$ we have $f_i(a_j)=\delta_{ij}$. This condition implies that  $f_1,f_2,\dots,f_r$ are linearly independent. Now, we show that these functions are contained in a subspace of dimension $\binom{n}{q-1}+\binom{n}{q-2}+\dots+\binom{n}{1}+\binom{n}{0}$.

Let $a_{i,j}$ denote the $j$th entry of $a_i$. 
Observe that for every $1\leq i\leq r$ and every $x\in A$
$$\langle x-a_i,x-a_i\rangle=\sum\limits_{j=1}^n (x_j-a_{i,j})^2=\sum\limits_{j=1}^n x_j^2-\sum\limits_{j=1}^n 2a_{i,j}x_j+\sum\limits_{j=1}^n a_{i,j}^2=2n-\sum\limits_{j=1}^n 2a_{i,j}x_j,$$
where the last equality holds, since $a_{i,j}\in \{\pm1\}$ and $x\in A\subseteq \{-1,1\}^n$.

Therefore, $f_i(x)=1-\left(2n-\sum\limits_{j=1}^n 2a_{i,j}x_j\right)^{q-1}$ is a polynomial of $x_1,x_2,\dots,x_n$ of degree $q-1$.

Furthermore, since $x_j^2=1$ for every $x_j\in A$ we can reduce the exponent of $x_j$ to 0 or 1, according to the parity of the original exponent. This way each $f_i$ is represented as a polynomial of degree at most $q-1$, where each individual degree is at most 1.

Each such monomial is uniquely determined with the subset of those variables $x_i$ which has exponent 1  (the rest of the variables have exponent 0). 
The number of these monomials is $$\binom{n}{q-1}+\binom{n}{q-2}+\dots+\binom{n}{1}+\binom{n}{0},$$
hence by the linear independency of $f_1,\dots,f_r$ this also serves as an upper bound for $r=|A|$, which completes the proof.

If $q\mid n$, then we can improve on this upper bound. Since, in this case $2n=0$ (in $\mathbb{F}_q$), thus 
$$f_i(x)=1-\left(2n-\sum\limits_{j=1}^n 2a_{i,j}x_j\right)^{q-1}=1-\left(\sum\limits_{j=1}^n 2a_{i,j}x_j\right)^{q-1}.$$
Expressing $f_i$ this way with the exception of the constant term 1 all the monomials have degree $q-1$. After reducing the exponents, we get monomials with even total degree (as $q-1$ is also even). Consequently, the upper bound is improved to the number of those subsets of $x_1,\dots,x_n$ whose size is even and at most $q-1$. Hence, $|A|=r\leq \binom{n}{q-1}+\binom{n}{q-3}+\dots+\binom{n}{2}+\binom{n}{0}$.

\end{proof}

\section{Concluding remarks} \label{sec-concl}

In this paper we prove that the largest subset of $\mathbb{F}_q^n$ (where $q$ is an odd prime) is between $\Theta(n^{q/3})$ and $\Theta(n^{q-2})$, however, for $q>3$ this leaves open what the right exponent is between $q/3$ and $q-2$. We shall note that one can transform our lower bound construction to a set $A$ for which $\langle x,x\rangle =0$ for every $x\in A$ and $\langle x,y\rangle \in B$ for any two distinct $x,y\in A$ with some sum-free subset $B\subseteq \mathbb{F}_q$. Then $\langle x-z,y-z\rangle =\langle x,y\rangle - (\langle x,z\rangle +\langle z,y\rangle)\ne 0$ indeed holds.  

For the analogous problem for $k$-right corners there's also a gap between our lower bound and Naslund's upper bound (which gap is getting larger as $k$ grows). In the case of triangles with all right angles, there's also a gap in the exponent of $n$.

For the problem where self-orthogonal differences are to be avoided our bounds are asymptotically tight, though the exact answer is determined only for a specific (infinite) family of parameters $q$ and $n$.

It would be interesting to further tighten these gaps.

\section{Acknowledgements}
The research was supported by the Lend\"ulet program of the Hungarian Academy of Sciences (MTA). MD was also supported by the New National Excellence Program of the National Research, Development and Innovation Fund and the Ministry for Innovation and Technology (\'UNKP-20-1). PPP was also supported by the National Research, Development and Innovation Office NKFIH (Grant Nr. K124171, K129335 and BME NC TKP2020). 

\medskip


\begin{thebibliography}{99}



\bibitem{Babai}
L.~Babai, H.~Snevily, R.~M.~Wilson,
A new proof of several inequalities on codes and sets,
Journal of Combinatorial Theory, Series A
71 (1) (1995) 146--153.




\bibitem{Bennett}
M.~Bennett,
Occurrence of right angles in vector spaces over finite fields,
European Journal of Combinatorics, {70} (2018) 155--163.


\bibitem{Croot-Lev-Pach}
E.~Croot, V.F.~Lev, P.P.~Pach,
Progression-free sets in $\Z_4^n$ are exponentially small, 
Ann. of Math. (2) 185 (2017), no. 1, 331--337. 

\bibitem{Delsarte}
P.~Delsarte, Four Fundamental Parameters of a Code and Their Combinatorial Significance, Information and Control 23 (1973) 407--438.
\bibitem{Ellenberg-Gijswijt}
J.~S.~Ellenberg, D.~Gijswijt, 
On large subsets of $\F_q^n$ with no three-ter arithmetic progression, Ann. of Math. (2) 185 (2017), no. 1, 339--343. 


\bibitem{EFF}
P.~Erd\H{o}s, P.~Frankl, Z.~F\"uredi,
Families of finite sets in which no set is covered by the union of two others,
Journal of Combinatorial Theory, Series A
33 (2) (1982) 158--166.


\bibitem{Frankl}
P.~Frankl,
Orthogonal vectors in the $n$-dimensional cube and codes with missing distances, Combinatorica 6 (3) (1986) 279--285.

\bibitem{Ge}
G.~Ge, C.~Shangguan,
{Maximum subsets of $\mathbb{F}_q^n$ containing no right angles},  J. Algebr. Comb. (2019). https://doi.org/10.1007/s10801-019-00908-4






\bibitem{Naslund}
E.~Naslund,
The partition rank of a tensor and $k$-right corners in $\mathbb{F}_q^n$,
Journal of Combinatorial Theory, Series A 174 (2020) 105190


\bibitem{Rodl}
V.~R\"odl,
On a Packing and Covering Problem,
European Journal of Combinatorics 6 (1) (1985) 69--78.



\bibitem{Tao2}
W.~Sawin and T.~Tao,
Notes on the slice rank of tensors,\\ {https://terrytao.wordpress.com/2016/08/24/notes-on-the-slice-rank-of-tensors/} (2016)

\bibitem{Tao1}
T.~Tao, A symmetric formulation of the Croot-Lev-Pach-Ellenberg-Gijswijt capset bound,\\ 
{https://terrytao.wordpress.com/2016/05/18/a-symmetric-formulation-of-the-croot-lev-pach-ellenberg-gijswijt-capset-bound} (2016)


\end{thebibliography}
\end{document}